# Some Counterexamples in Ring Theory

Santosh Kumar Pandey,

Faculty of Science, Technology and Forensic,

SPUP, Opp. Vigyan Nagar, Jodhpur-342037, India.

Email: skpandey12@gmail.com

**ABSTRACT**

In this note we provide important and significant observations in ring theory related to weakly tripotent rings. We provide counterexamples for the structure theorem for commutative weakly tripotent rings appeared in arXiv (2017) and Bull. Korean Math. Soc.(2018). Also we give some additional results to improve the existing results.

## 1. Introduction

Let $R$ be an associative ring with identity. Let $J(R)$ denotes the Jacobson radical of $R$ and $U(R)$ stands for the set of units in $R$. A ring $R$ satisfying $a^3 = a$ for each $a \in R$ is widely known as tripotent ring [1-3]. Similarly a ring $R$ is called weakly tripotent if $a^3 = a$ or $(1+a)^3 = 1+a$ for each $a \in R$ [4-5].

Evidently it follows that each tripotent ring is a weakly tripotent ring but a weakly tripotent ring need not be a tripotent ring.

In the mathematical literature one may find the following important structure theorem for commutative weakly tripotent rings [4-5].

A commutative ring $R$ is a weakly tripotent ring iff $R = R_1 \times R_2$ such that

(1) $R_2$ is a tripotent ring of characteristic three or $R_2 = 0$;

(2) $R_1 = 0$ or $3 \in U(R_1)$ and $R_1$ can be embedded as a subring of a direct product $R_0 \times (\prod_{i \in I} R_i)$ such that $R_0$ is a weakly tripotent ring without nontrivial idempotents, and all $R_i$ are Boolean rings.

It may be emphasized that as per [4-5], for each commutative weakly tripotent ring $R$ there exists a non-zero ideal $L$ maximal with the property $L \cap J(R) = 0$ and it is necessary for $R_1$ to be embedded as a subring of a direct product $R_0 \times (\prod_{i \in I} R_i)$ such that $R_0$ is a weakly tripotent ring without nontrivial idempotents, and all $R_i$ are Boolean rings.

However in this note we exhibit the following observations:

There exist rings $R$ having the following properties: (a) $R$ is a commutative weakly tripotent ring. (b) $R$ does not have any non-zero ideal $L$ maximal with the property $L \cap J(R) = 0$ (c) $R$ is decomposable as $R = R_1 \times R_2$, with $R_2 = 0$, $3 \in U(R_1)$ and $R_1$ is embedded as a subring of a direct product $R_0 \times (\prod_{i \in I} R_i)$ such that $R_0$ is a weakly tripotent ring without nontrivial idempotents, and all $R_i$ are Boolean rings. In addition we give a rectification of this discrepancy which was unnoticed since 2017.

Following [4-5] in this note we shall consider commutative weakly tripotent ring of characteristic $2^t$, $t \in \{1, 2, 3\}$.

## 2. Observations

Let $R$ be a commutative weakly tripotent ring. As per [4-5] if we consider the set of all ideals of $R$, then this set will contain a maximal element $L$ such that $L \cap J(R) = 0$ and for $0 \neq a \in L$ we can find an ideal $I_a$ . Here $I_a$ is maximal with the property $a \notin I_a$ and by Zorn's Lemma $I_a$ can be chosen such that $J(R) \subseteq I_a$. Clearly as per [4-5] $L$ is a non-zero ideal. However we exhibit that there are rings $R$ of characteristic two, four and eight also such that these facts do not hold in them even though they are weakly tripotent having the decomposition as described below.

More explicitly in the rings described in the following Theorem 2.1, Theorem 2.2 and Theorem 2.3 there does not exist any non-zero ideal $L$ maximal with the property $L \cap J(R) = 0$ and so $I_a$ doest not exist such that $J(R) \subseteq I_a$ with $a \notin I_a$, where $0 \neq a \in L$. Nontheless $R$ is decomposable as $R = R_1 \times R_2$, with $R_2 = 0$, $3 \in U(R_1)$ and $R_1$ is embedded as a subring of a direct product $R_0 \times (\prod_{i \in I} R_i)$ such that $R_0$ is a weakly tripotent ring without nontrivial idempotents, and all $R_i$ are Boolean rings.

**Theorem 2.1:** There exists a ring $R$ of characteristic $2$ having the following properties: (a) $R$ is a commutative weakly tripotent ring. (b) $R$ does not have any non-zero ideal $L$ maximal with the property $L \cap J(R) = 0$ (c) $R$ is decomposable as $R = R_1 \times R_2$, with $R_2 = 0$, $3 \in U(R_1)$ and $R_1$ is embedded as a subring of a direct product $R_0 \times (\prod_{i \in I} R_i)$ such that $R_0$ is a weakly tripotent ring without nontrivial idempotents, and all $R_i$ are Boolean rings.

**Proof.** Let $R = R_1 \times R_2$. Here $R_2 = 0$ and $R_1 = \{0, a, b, c\}$ is the commutative ring of order four and characteristic two having only one non-zero nilpotent element. It is well known that there is only one such a commutative ring of order four and characteristic two. We have $a + b = c, a + c = b, b + c = a$ and $a^2 = b^2 = a, c^2 = 0$, $ab = b, ac = c, bc = c$, etc.

Clearly $R = R_1 \times R_2 = \{(0,0), (a,0), (b,0), (c,0)\}$ is a commutative ring of characteristic two. One may verify that

$$(0,0)^3 = (0,0),\ (a,0)^3 = (a,0), (b,0)^3 = (b,0) \text{ and } [(a,0) + (c,0)]^3 = (a,0) + (c,0).$$

Therefore $R$ is a commutative weakly tripotent ring. One may further note that the Jacobson radical $J(R)$ of $R$ coincides with the nil radical of $R$. Also, one may easily find that $R$ does not have any non-zero ideal $L$ maximal with the property $L \cap J(R) = 0$ and hence $I_a$ doest not exist such that $J(R) \subseteq I_a$ with $a \notin I_a$ where $0 \neq a \in L$. Even though $R_1$ is embedded as a subring of a direct product $R_0 \times (\prod_{i \in I} R_i)$ such that $R_0$ is a weakly tripotent ring without nontrivial idempotents, and all $R_i$ are Boolean rings. Thus the proof is complete.

**Theorem 2.2:** There exists a ring $R$ of characteristic $4$ having the following properties: (a) $R$ is a commutative weakly tripotent ring. (b) $R$ does not have any non-zero ideal $L$ maximal with the property $L \cap J(R) = 0$ (c) $R$ is decomposable as $R = R_1 \times R_2$, with $R_2 = 0$, $3 \in U(R_1)$ and $R_1$ is embedded as a subring of a direct product $R_0 \times (\prod_{i \in I} R_i)$ such that $R_0$ is a weakly tripotent ring without nontrivial idempotents, and all $R_i$ are Boolean rings.

**Proof.** Let $R = R_1 \times R_2$. Here $R_2 = 0$ and $R_1 = Z_4$. Clearly $R = \{(0,0),(1,0),(2,0),(3,0)\}$ is a commutative ring of characteristic four.

We note that $(0,0)^3 = (0,0)$, $(1,0)^3 = (1,0)$, $(3,0)^3 = (3,0)$ and $[(1,0)+(2,0)]^3 = (3,0)$ $= (1,0)+(2,0)$. Therefore $R$ is a commutative weakly tripotent ring.

One may further note that the Jacobson radical $J(R)$ of $R$ coincides with the nil radical of $R$. One may easily find that $R$ does not have any non-zero ideal $L$ maximal with the property $L \cap J(R) = 0$ and hence $I_a$ doest not exist such that $J(R) \subseteq I_a$ with $a \notin I_a$ where $0 \neq a \in L$. However $R_1$ is embedded as a subring of a direct product $R_0 \times (\prod_{i \in I} R_i)$ such that $R_0$ is a weakly tripotent ring without nontrivial idempotents, and all $R_i$ are Boolean rings. Hence the proof is complete.

**Theorem 2.3:** There exists a ring $R$ of characteristic $8$ having the following properties: (a) $R$ is a commutative weakly tripotent ring. (b) $R$ does not have any non-zero ideal $L$ maximal with the property $L \cap J(R) = 0$ (c) $R$ is decomposable as $R = R_1 \times R_2$, with $R_2 = 0$, $3 \in U(R_1)$ and $R_1$ is embedded as a subring of a direct product $R_0 \times (\prod_{i \in I} R_i)$ such that $R_0$ is a weakly tripotent ring without nontrivial idempotents, and all $R_i$ are Boolean rings.

**Proof.** Let $R = R_1 \times R_2$. Here $R_2 = 0$ and $R_1 = Z_8$. We leave it for readers to verify that $R$ is a commutative weakly tripotent ring. Also we note that $R$ does not have any non-zero ideal $L$ maximal with the property $L \cap J(R) = 0$ and hence $I_a$ doest not exist such that $J(R) \subseteq I_a$ with $a \notin I_a$ where $0 \neq a \in L$. Nevertheless $R_1$ is embedded as a subring of a direct product $R_0 \times (\prod_{i \in I} R_i)$ such that $R_0$ is a weakly tripotent ring without nontrivial idempotents, and all $R_i$ are Boolean rings. Thus the result follows.

Below we shall present rectification and discussion.

## 3. **Rectification and Discussion**

Now we shall rectify the above issue which was unnoticed since 2017 [4-5]. Let $R$ be a commutative weakly tripotent ring of characteristic $2^t$, $t \in \{1,2,3\}$. We have seen that if $R$ does not have any non-zero ideal $L$ maximal with the property $L \cap J(R) = 0$, then $I_a$ doest not exist such that $J(R) \subseteq I_a$ with $a \notin I_a$, where $0 \neq a \in L$. In this case $J(R)$ is the unique maximal ideal of $R$ and $R$ does not have any non-zero ideal $L$ such that $L \cap J(R) = 0$. Now we have the following result for commutative weakly tripotent rings in which $J(R)$ is the unique maximal ideal.

**Theorem 3.1:** If $R$ is a commutative weakly tripotent ring in which $J(R)$ is the unique maximal ideal then $R \cong R \times \{0\} = A \times B$ where $A$ is a weakly tripotent ring without non-trivial idempotents and $B$ is a Boolean ring of characteristic one.

**Proof.** Let $R$ is a commutative weakly tripotent ring in which $J(R)$ is the unique maximal ideal. This implies that $R$ is a commutative indecomposable ring. It may be noted that an indecomposable commutative ring is a ring without non-trivial idempotents. Hence $R \cong R \times \{0\} = A \times B$, where $A$ is a weakly tripotent ring without non-trivial idempotents and $B$ is a Boolean ring of characteristic one.

**Theorem 3.2:** If $R$ is a commutative ring in which $J(R)$ is the unique maximal ideal such that $R \cong R \times \{0\} = A \times B$ where $A$ is a ring without non-trivial idempotents and $B$ is a Boolean ring of characteristic one. Then $R$ need not be a weakly tripotent ring.

Proof. . This can be seen as follows.

Let $R = \left\{ \begin{pmatrix} 0 & 0 \\ 0 & 0 \end{pmatrix}, \begin{pmatrix} 1 & 0 \\ 0 & 1 \end{pmatrix}, \begin{pmatrix} 0 & 1 \\ 1 & 1 \end{pmatrix}, \begin{pmatrix} 1 & 1 \\ 1 & 0 \end{pmatrix} \right\}$. One may note that $R$ is a ring under addition and multiplication of matrices modulo two and $J(R)$ is the unique maximal ideal of $R$. Also $R \cong R \times \{0\} = A \times B$ where $A$ is a commutative ring without non-trivial idempotents and $B$ is a Boolean ring of characteristic one. However $R$ is not a weakly tripotent ring. Thus the proof is complete.

Hence if $R$ is a commutative weakly tripotent ring in which $J(R)$ is the unique maximal ideal then there does not exist a non-zero ideal $L$ maximal with the property $L \cap J(R) = 0$ and $I_a$ with $a \notin I_a$, where $0 \neq a \in L$ such that $J(R) \subseteq I_a$. Even though $R$ is decomposable as $A \times B$ .where $A$ is a weakly tripotent ring without non-trivial idempotents and $B$ is a Boolean ring of characteristic one.

**Acknowledgements.** The author is thankful to A. Pandit for his kind support.

**Statement and Declaration:** The author declares that there is no competing interest and no funding support was available for this research.